\documentclass[12pt]{amsart}
\usepackage{psfrag}
\usepackage{mathrsfs}
\usepackage{amssymb}

 \usepackage[notcite,notref,color,final]{showkeys}
 \definecolor{refkey}{rgb}{0,0,1}
 \definecolor{labelkey}{rgb}{0,0,1}
\usepackage{mathpazo}

\usepackage{hyperref}
\usepackage{graphicx, xspace, color}

\theoremstyle{plain}
\newtheorem{lemma}{Lemma}
\newtheorem{theorem}[lemma]{Theorem}

\newtheorem{corollary}[lemma]{Corollary}
\newtheorem{example}[lemma]{Example}

\newcommand{\etal}{\emph{et al.}\xspace}
\newcommand{\shift}[2]{{#1}^{\langle{#2}\rangle}}
\newcommand{\jsh}[3]{#2^{\langle{#1}\rangle}}

\newcommand{\bdry}[1]{\partial{#1}}

\newcommand{\seq}[2]{\ensuremath{#1_0, \ldots, #1_{#2-1}}}

\newcommand{\setbp}{\mathscr{B}}

\newcommand{\setap}[1]{\mathscr{U}^{#1}}
\newcommand{\setau}{\mathscr{U}}

\newcommand{\ba}{\ensuremath{ {\bf a}}}

\newcommand{\pth}[1]{\mathcal{#1}}
\newcommand{\setbb}[2]{\mathscr{B}(#2, #1)}
\newcommand{\setgg}[2]{\mathscr{G}(#2, #1)}
\newcommand{\setaa}[2]{\mathscr{A}(#2, #1)}

\newcommand{\bee}{b}
\newcommand{\Bee}{B}

\author{J. Irving}
\address{Department of Mathematics and Computing Science \\
        St. Mary's University, Halifax, NS, B3H 3C3, Canada}
\email{john.irving@smu.ca}

\author{A. Rattan}
\address{Department of Mathematics, Room 2-376, Massachusetts Institute of Technology,
Cambridge, MA, 02139-4307, USA}
\email{arattan@math.mit.edu}

\title[Cyclically shifting boundary]{The number of lattice
paths below a cyclically shifting boundary}

\begin{document}

\begin{abstract}
    We count the number of lattice paths lying under a cyclically
    shifting piecewise linear boundary of varying slope.  Our main result
    extends well known enumerative formulae concerning lattice paths, and its derivation involves a classical
    reflection argument. A refinement allows for the counting of paths with a specified number of corners.
    We also apply the result to examine paths dominated  by periodic boundaries.
\end{abstract}

\maketitle
\section{Introduction}
\label{sec:introduction}


Throughout, the term \emph{lattice path} refers to a path in the integer lattice $\mathbb{Z} \times \mathbb{Z}$
with unit steps up and to the right (\emph{i.e.} steps $(0,1)$ and $(1,0)$, respectively).

\begin{figure}
    \includegraphics[width=.35\textwidth]{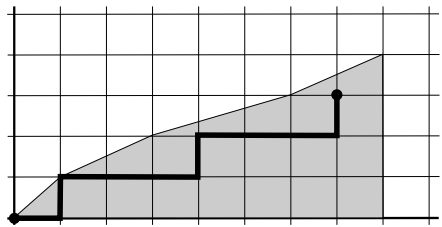}
    \caption{A path dominated by $\ba=(1,2,3,2)$.}\label{fig:pathexmp}
\end{figure}

Let $\ba = (a_0,\ldots,a_{m-1})$ be a weak $m$-part composition of $n$ (recall that this means the $a_i$ are
nonnegative integers summing to $n$). This paper concerns the enumeration of lattice paths from the origin that
lie weakly under the piecewise linear boundary curve $\bdry{\ba}$ defined by
$$
    x =
    a_i y + \sum_{j=0}^{i-1} a_j,
    \quad\text{for $y \in [i,i+1]$}.
$$
We say such paths (and all points weakly under $\bdry{\ba}$) are \emph{dominated} by $\ba$.  The boundary
corresponding to $\ba=(1,2,3,2)$ is shown in Figure~\ref{fig:pathexmp}, along with a path it dominates.

\newcommand{\domset}[1]{D(#1)}

\begin{figure}
    \includegraphics[width=.8\textwidth]{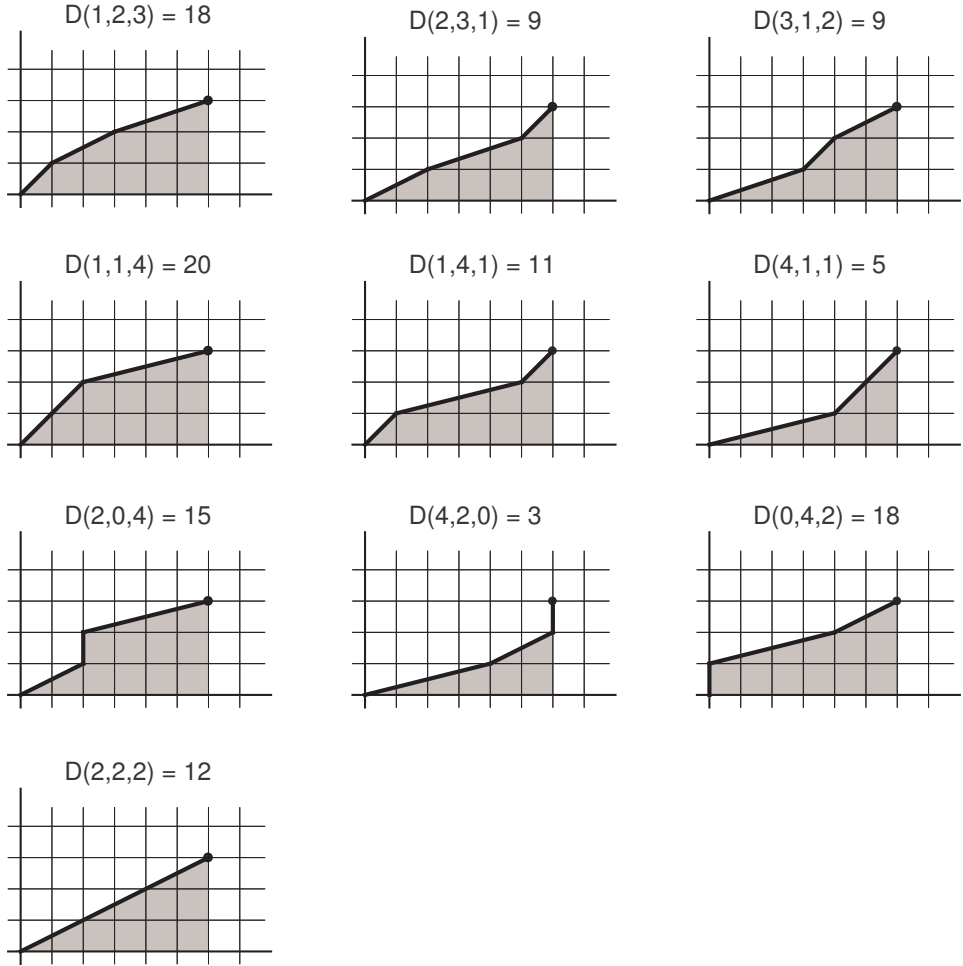}
    \caption{The number of paths dominated by cyclically shifted boundaries.}\label{fig:exampboundaries}
\end{figure}

Let $\domset{\ba}$ be the number of paths from $(0,0)$ to $(n,m)$ dominated by $\ba$.  For example, the numbers
$\domset{\ba}$ for various 3-part compositions of 6 are given above their boundaries in
Figure~\ref{fig:exampboundaries}.
When all parts of $\ba$ are the same it is well-known~\cite[Exercise 5.3.5]{goul:1} that $D(\ba)$ is a generalized
Catalan number: We have
\begin{align}
    \domset{\underbrace{a,a,\ldots,a}_{\text{$m$ copies}}} =
        \frac{1}{(a+1)m+1} \binom{(a+1)m+1}{m},
\end{align}
where the case $a=1$  corresponds with the classical Dyck paths counted by the usual Catalan numbers.  However,
for general $\ba$ no simple formula for $\domset{\ba}$ is known, and indeed it is unlikely that such a formula
exists (though the Kreweras dominance theorem~\cite[Section 5.4.7]{goul:1} does give a determinantal expression).
It is the purpose of this paper to show that simple enumerative formulae do hold provided that we consider paths
dominated by \emph{all cyclic shifts} of an arbitrary composition.

For an integer $j$, let $\shift{\ba}{j}$ denote the $j$-th cyclic shift of $\ba$, namely
\begin{equation}
\label{eq:shift}
    \shift{\ba}{j} = (a_{-j},a_{-j+1},\ldots,a_{-j+m-1}),
\end{equation}
where the indices are to be interpreted modulo $m$.  For example, the rows of Figure~\ref{fig:exampboundaries}
show $\bdry{\ba}, \bdry{\shift{\ba}{1}}$, and $\bdry{\shift{\ba}{2}}$ for the compositions $\ba=(1,2,3), (1,1,4),
(2,0,4)$ and $(2,2,2)$.  Notice that in each of the top three rows there are a total of 36 dominated paths from
$(0,0)$ to $(6,3)$, and this many also in the bottom row if the three identical cyclic shifts of $\ba=(2,2,2)$ are
taken into account. That is, $D(\ba)+D(\shift{\ba}{1})+D(\shift{\ba}{2})=36$ for each of these 3-part compositions
$\ba$ of 6.  This is a special case of a more general phenomenon.

We define a \emph{lattice path boundary pair (LPBP)} to be an ordered pair $(\pth{P}, (\ba,j))$, where $\pth{P}$
is a lattice path beginning at the origin, $\ba$ is a weak $m$-part composition, and $j$ is an integer with $0
\leq j < m$. If $\pth{P}$ is dominated by $\shift{\ba}{j}$ then we say $(\pth{P}, (\ba,j))$ is a \emph{good pair},
otherwise it is a \emph{bad pair}. Let $\setaa{t}{\ba}$ be the set of all LPBPs of the form $(\pth{P}, (\ba,j))$,
where $\pth{P}$ terminates at the point $t$. Let $\setbb{t}{\ba}$ and $\setgg{t}{\ba}$ be the subsets of
$\setaa{t}{\ba}$ consisting of bad and good pairs, respectively. Clearly, $\setaa{t}{\ba} = \setbb{t}{\ba} \cup
\setgg{t}{\ba}$, with the union  disjoint.

We are now ready to state our main result.  After its discovery, we found an essentially equivalent conjecture in
earlier work of Tamm~\cite{tamm}. Though Tamm's paper concerns paths under periodic boundaries (see
Section~\ref{sec:periodic}), the conjecture itself is coarsely formulated in the language of \emph{two-dimensional
arrays}, with a proof only in the case $m=2$.

\begin{theorem}\label{thm:gen}
    Let $\ba$ be a weak $m$-part composition of $n$ and let $t=(k,l)$, with
    $0 \leq k \leq n$,  $0 \leq l \leq m$. If the point $(k+1,l)$ lies weakly to the right of $\bdry{\ba}$ then
    \begin{equation}\label{eq:trivial}
        |\setaa{t}{\ba}| = m {k + l \choose l},
    \end{equation}
    \begin{equation}\label{eq:thm1}
        |\setbb{t}{\ba}|  = n\binom{k+l}{l-1},
    \end{equation}
    and
    \begin{equation}\label{eq:breakdown}
        |\setgg{t}{\ba}|  = |\setaa{t}{\ba}| - |\setbb{t}{\ba}| =
        \frac{m(k+1)-nl}{k+1}\binom{k+l}{l}.
    \end{equation}
    \qed
\end{theorem}

Thus we have the surprising fact that the total number of paths dominated by all cyclic shifts of a piecewise
linear boundary does not depend on the specific parts of its defining composition $\ba$.  Instead, allowing all
shifts the boundary acts as an averaging process with a very pleasant enumerative outcome.

Clearly the hypothesis of Theorem~\ref{thm:gen} are satisfied by any terminus $(k,l)$ that is itself dominated by
all cyclic shifts of $\ba$.  In particular, setting $(k,l)=(n,m)$ in the theorem explains our previous observation
that there are $36 = \binom{6+3}{3-1}$ dominated paths for each row of Figure~\ref{fig:exampboundaries}.

\begin{corollary}\label{cor:total}  For any weak $m$-part composition $\ba$ of $n$, we have
$$
\domset{\ba} +\domset{\shift{\ba}{1}}+\cdots+\domset{\shift{\ba}{m-1}}= \binom{n+m}{m-1}.
$$
\qed
\end{corollary}

Now consider the composition $\ba=(a,a,\ldots,a)$ of $n=ma$. Observe that $\shift{\ba}{i}=\ba$ for all $i$, while
$\bdry{\ba}$ is simply the line $x=ay$. Applying Theorem~\ref{thm:gen} and dividing by $m$ to remove the effect of
boundary rotation therefore yields the following well-known result, often referred to as the \emph{generalized
ballot theorem}. See the survey article~\cite{kratt:1} for more information.
\begin{corollary}
\label{cor:kratt} If $k \geq al$, then there are
\begin{equation*}\label{eq:kratt}
    \frac{k- al + 1}{k+1}{k+l \choose l}
\end{equation*}
lattice paths from $(0,0)$ to $(k,l)$ that lie weakly below the line $x=ay$. \qed
\end{corollary}

In the next section we give a bijective proof of Theorem~\ref{thm:gen}. Of course, this amounts to
proving~\eqref{eq:thm1} since~\eqref{eq:trivial} is trivial. We prove~\eqref{eq:thm1} by showing bad paths are in
bijection with a less restrictive set of paths, in the spirit of Andr\'e's~\cite{andre} reflection principle. In
fact, our proof is a generalization of the bijection used in~\cite{goulserr:1} to prove Corollary~\ref{cor:kratt}.
(Our bijection reduces to that of~\cite{goulserr:1} in the case when all parts of $\ba$ are the same, though we
must make an allowance for the cyclically shifting boundary.)

Section~\ref{sec:cyclelemma} contains a brief account of an alternative derivation of Theorem~\ref{thm:gen} using
the Cycle Lemma.   In Section~\ref{sec:refinements} we present a refinement of the theorem that counts paths with
a specified number of corners. Finally, Section~\ref{sec:periodic} illustrates a handful of applications of the
theorem to the enumeration of lattice paths lying under periodic boundaries. Interestingly, all of these
applications pivot on the fact that the hypotheses of Theorem~\ref{thm:gen} do not require the terminus $(k,l)$ to
be dominated by all cyclic shifts of $\ba$.

\section{A Proof of Theorem \ref{thm:gen}}
\label{sec:proof}

\newcommand{\wlwb}{\leq}
\newcommand{\wlb}{\lesssim}
\newcommand{\lb}{<}


Throughout this section we have in mind a fixed weak composition $\ba = (\seq{a}{m})$ of $n$ and its corresponding
boundary $\bdry{\ba}$.  For arbitrary $j \in \mathbb{Z}$ we interpret the symbol $a_j$ to mean $a_{j\!\! \mod m}$.

\begin{figure}
    \psfrag{pst:S}{$\jsh{2}{p}{\ba}$}
    \psfrag{pst:p}{$p$}
    \psfrag{pst:a}{$\bdry{\ba}$}
    \psfrag{pst:as}{$\bdry{\shift{\ba}{2}}$}
    \includegraphics[width=.4\textwidth]{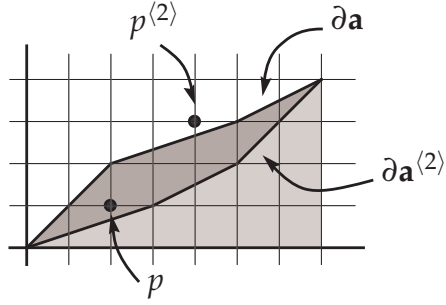}
    \caption{A point $p$ and its shift $\jsh{2}{p}{\ba}$.}\label{fig:shifting}
\end{figure}

For any lattice point $p=(x,y)$ with $0 \leq x < n$ and $1 \leq y \leq m$, and for any integer $j$,  define
the \emph{$j$-th shift of $p$ (relative to $\ba$)} to be the point
\begin{equation}
    \label{eq:pshift}
    \jsh{j}{p}{\ba} = (x + a_{-1}+a_{-2}+\cdots+a_{-j}\!\! \mod n,~~ y + j\!\! \mod m),
\end{equation}
where the modular reductions in the first and second coordinate are understood to yield representatives in
$\{0,1,\ldots,n-1\}$ and $\{1,2,\ldots,m\}$, respectively.  Informally, $\jsh{j}{p}{\ba}$ is in the same position
relative to $\bdry{\shift{\ba}{j}}$ as $p$ is to $\bdry{\ba}$. (See Figure~\ref{fig:shifting}.)

Finally, we define the relations $\wlwb$, $\wlb$, and $\lb$ on  lattice points as follows:
\begin{itemize}
\item   $(x_1,y_1) \wlwb (x_2,y_2) \Leftrightarrow x_1 \leq x_2$ and $y_1 \leq y_2$
\item   $(x_1,y_1) \wlb (x_2,y_2) \Leftrightarrow x_1 \leq x_2$ and $y_1 < y_2$
\item   $(x_1,y_1) \lb (x_2,y_2) \Leftrightarrow x_1 < x_2$ and $y_1 < y_2$
\end{itemize}
\vspace*{3mm}

\begin{figure}
    \psfrag{ps:pi}{$p_1$}
    \psfrag{ps:pi1}{$\jsh{1}{p_1}{\ba}$}
    \psfrag{ps:pi2}{$\jsh{2}{p_1}{\ba}$}
    \psfrag{ps:pi3}{$\jsh{3}{p_1}{\ba}$}
    \psfrag{ps:ti}{$p_4$}
    \psfrag{ps:ti1}{$\jsh{1}{p_4}{\ba}$}
    \psfrag{ps:ti2}{$\jsh{2}{p_4}{\ba}$}
    \psfrag{ps:ti3}{$\jsh{3}{p_4}{\ba}$}
    \includegraphics[width=.8\textwidth]{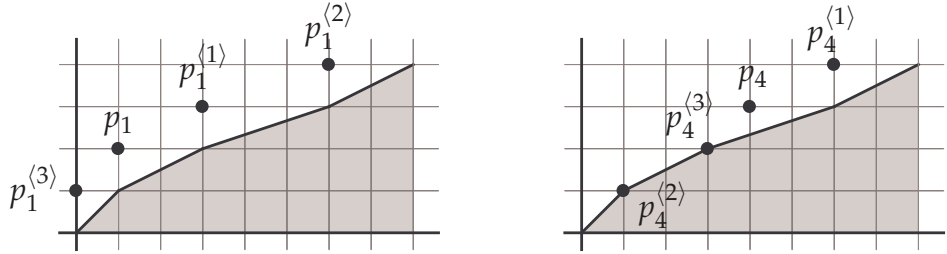}
    \caption{The sets $\Bee_1$ and $\Bee_4$ relative to the composition $\ba=(1,2,3,2)$.}\label{fig:bsets}
\end{figure}

For $0 \leq i < n$, let $p_i=(i,y_i)$, where $y_i$ is the least
integer such that the point $(i,y_i)$ lies strictly above $\bdry{\ba}$.  Define
\begin{equation}\label{eq:setbi}
    \Bee_i := \{ p_i, \jsh{1}{p_i}{\ba}, \jsh{2}{p_i}{\ba}, \ldots, \jsh{m-1}{p_i}{\ba} \}.
\end{equation}
to be the set of lattice points having the same relative positions to the boundary curves
$\bdry{\ba},\bdry{\shift{\ba}{1}},\ldots,\bdry{\shift{\ba}{m-1}}$ as the point $p_i$ has to $\ba$. (See
Figure~\ref{fig:bsets}.)  Note that the sets $\Bee_0,\ldots,\Bee_{m-1}$ are not disjoint.

Let $\setbp_i(\ba,t)$ be the set of all bad LPBPs of the form $(\pth{P},(\ba,j))$, where the path $\pth{P}$
terminates at $t$ and its first \emph{bad} step (\emph{i.e} the first step crossing $\bdry{\ba^j}$) lands at
the point $\jsh{j}{p_i}{\ba} \in \Bee_i$.  Then clearly $\setbb{t}{\ba} = \bigcup_{i=0}^{n-1}
\setbp_i(\ba,t)$, with the union being disjoint.  We shall prove Theorem~\ref{thm:gen} by showing that
$|\setbp_i(\ba,t)|$ is independent of $i$.

\begin{figure}
\begin{center}
    \psfrag{ps:pi}{$\bee_1^1$}
    \psfrag{ps:pi1}{$\bee_1^2$}
    \psfrag{ps:pi2}{$\bee_1^3$}
    \psfrag{ps:pi3}{$\bee_1^0$}
    \psfrag{ps:ti}{$\bee_4^2$}
    \psfrag{ps:ti1}{$\bee_4^3$}
    \psfrag{ps:ti2}{$\bee_4^0$}
    \psfrag{ps:ti3}{$\bee_4^1$}
    \includegraphics[width=.75\textwidth]{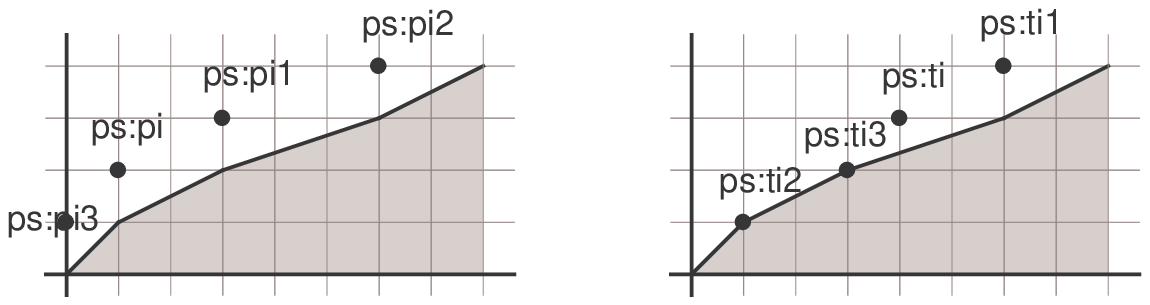}
    \caption{Construction of $\Bee_1=\{\bee_1^j\}$ and $\Bee_4=\{\bee_4^j\}$ relative to $\ba=(1,2,3,2)$.}\label{fig:bji}
\end{center}
\end{figure}

Observe that no two points in any given set $\Bee_i$ can have the same $y$ coordinates. In fact, let
$s_i=m+1-y_i$, so that $\jsh{s_i}{p_i}{\ba}$ has $y$-coordinate 1, and define
\begin{equation}
    \label{eq:bdefn}
    \bee^j_i = \jsh{s_i+j}{p_i}{\ba}, \qquad \text{for $0 \leq j < m$.}
\end{equation}
Then
\begin{equation}
    \Bee_i = \{\bee_i^0, \bee_i^1, \ldots, \bee_i^{m-1}\},
\end{equation}
where the $y$-coordinate of $\bee_i^j$ is $j+1$ and $\bee_i^0 \wlb \bee_i^{1} \wlb \cdots \wlb \bee_i^{m-1}$.
For example, for the sets $\Bee_1$ and $\Bee_4$ of Figure~\ref{fig:bsets} we have $s_1 =3$ and $s_4=2$,
respectively, and the appropriate relabellings  are shown in Figure~\ref{fig:bji}.

Let $0 \leq k \leq n$ and $0 \leq l \leq m$. We say $\Bee_i$ is \emph{complete with respect to the point
$t=(k,l)$} if  $\bee^{l-1}_i$ is weakly to the left of $t$.    For instance, the set $\Bee_1$ in
Figure~\ref{fig:bji} is complete with respect to $t=(3,3)$, while $\Bee_4$ is not.  The motivation behind this
definition will be made clear in the proof of the following lemma.

\begin{lemma}
    \label{thm:countbadpaths}
    If $\Bee_i$ is complete with respect to $t=(k,l)$ then
    \begin{equation*}
        |\setbp_i(\ba,t)| = {k + l \choose l-1}.
    \end{equation*}
\end{lemma}

\noindent \textbf{Proof.} We shall give a reflection-type correspondence between $\setbp_i(\ba,t)$ and the
set $\setau$ of all lattice paths from $(-1, 1)$ to $t$.

Let $\setbp_i^j \subseteq \setbp_i(\ba,t)$ consist of those LPBPs in which the first bad step lands at the
point $\bee_i^j$.  Since the $y$-coordinate of $\bee_i^j$ is $j+1$, clearly $\setbp_i^j = \emptyset$ for $j
\geq l$. Thus $\setbp_i(\ba,t) = \bigcup_{j=0}^{l-1} \setbp_i^j$, with the union disjoint.

Since $\Bee_i$ is complete with respect to $t$, we have $\bee^0_i \wlb \cdots \wlb \bee^{l-1}_i \wlwb t$.
That is, on each of the lines $y=1,\ldots,y=l$ there is a point in $\Bee_i$ that is weakly to the left of
$t$. It follows that any path from $(-1,1)$ to $t$ must intersect one of these points.  For $0 \leq j < l$,
let $\setau^j$ be the set of paths from $(-1,1)$ to $t$ that avoid the points
$\bee_i^0,\bee_i^1,\ldots,\bee_i^{j-1}$ but meet $\bee_i^j$.  Then, by our previous comment, $\setau =
\bigcup_{j=0}^{l-1} \setau^j$, with the union disjoint.

\begin{figure}
    \includegraphics[width=.9\textwidth]{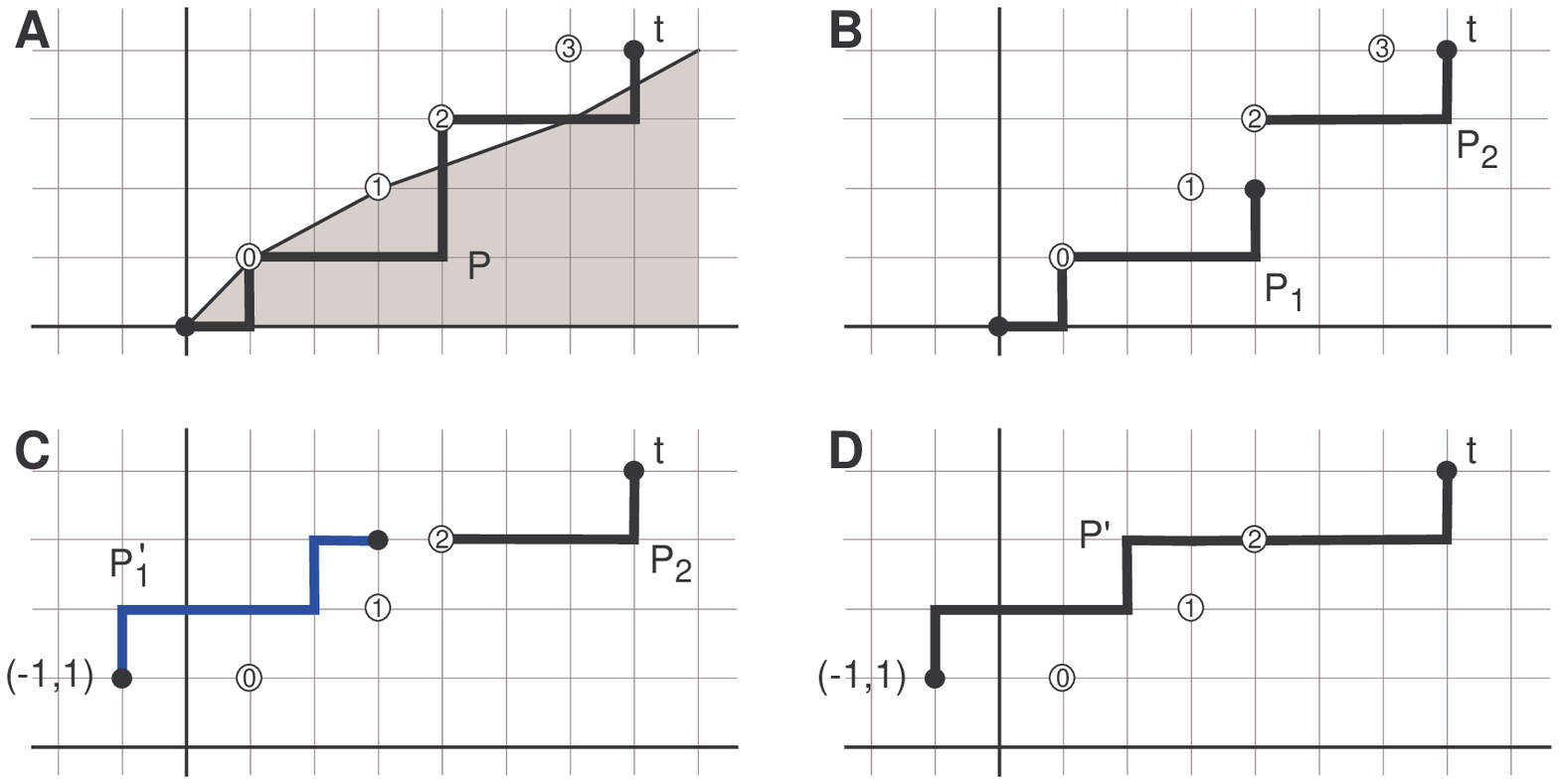}
    \caption{An illustration of the map $\psi_j : \setbp_i^j \longrightarrow \setap{j}$. Here $\ba=(1,2,3,2)$,
    $i=4$, $j=2$, and $s_i = 2$. The boundary
    $\bdry{\shift{\ba}{s_i+j}}=\bdry{\shift{\ba}{4}}=\bdry{\ba}$ is shown in panel
    A. Throughout, the point $b_i^j$ is indicated with a circled $j$.} \label{fig:bijection}
\end{figure}

We now define a mapping $\psi_j : \setbp_i^j \longrightarrow \setap{j}$ for each $j=0,\ldots,l-1$.  Given $L
\in \setbp_i^j$, construct $\psi_j(L)$ as follows. (See Figure~\ref{fig:bijection} for an illustration of the
construction.)
\begin{itemize}
\item[A.]   We have $\bee_i^j = (x,j+1)$ for some $x$, and $L=(\pth{P},(\ba,s_i+j))$ for some path
$\pth{P}$ whose first bad step is an up-step from $(x,j)$ to $\bee_i^j$.

\item[B.]   Remove this step to break $\pth{P}$ into two parts: the first, $\pth{P}_1$, is a path from $(0,0)$ to
$(x,j)$, and the second, $\pth{P}_2$, is a path from $(x,j+1)$ to $(k,l)$.

\item[C.]   Rotate $\pth{P}_1$ through $180^\circ$ and translate to obtain a new path $\pth{P}_1'$ beginning at $(-1, 1)$
and terminating at $(x-1,j+1)$.

\item[D.]   Join $\pth{P}_1'$ and $\pth{P}_2$ by adding a right-step from $(x-1,j+1)$ to $\bee_i^j$, thus creating a path
$\pth{P}'$ from $(-1,1)$ to $(k,l)$. Finally, set $\psi_j(L) = \pth{P}'$.
\end{itemize}

To ensure $\psi_j$ is well defined we must check that indeed $\pth{P}' \in \setap{j}$.  The only contentious
issue here is whether $\pth{P}'$ avoids the points $\bee_i^0,\ldots,\bee_i^{j-1}$. To see why this is the
case, consider the piecewise linear curve $\mathcal{C}$ obtained by joining the points $\bee_i^{-1},\bee_i^0,
\ldots, \bee_i^j$, where
\begin{equation}
\label{eq:bn}
    \bee_i^{-1} := ( - a_{-s_i} + (i+ a_{-1} + a_{-2} + \cdots + a_{-s_i} \!\!\! \mod n) ,\,\, 0).
\end{equation}
\begin{figure}
\begin{center}
    \includegraphics[width=\textwidth]{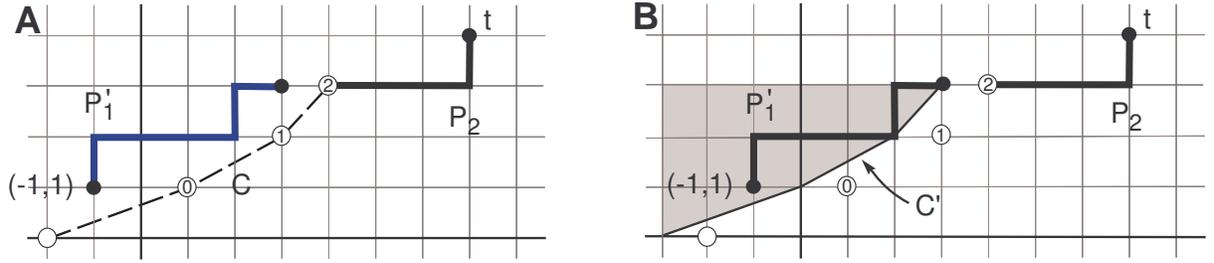}
    \caption{Proving that the map $\psi_2 : \setbp_4^2 \longrightarrow \setap{2}$
    illustrated in Figure~\ref{fig:bijection} is well defined.  The point $b_i^{-1}$ is
    indicated with an open circle.  The shaded region in panel B fits perfectly into
    that of Figure~\ref{fig:bijection}A after a $180^\circ$ rotation.}\label{fig:bijproof}
\end{center}
\end{figure}
See Figure~\ref{fig:bijproof}A for an illustration.

Since~\eqref{eq:pshift} and~\eqref{eq:bdefn} give
$$
    \bee_i^r = (i + a_{-1} + a_{-2} + \cdots + a_{-(s_i+r)} \!\!\mod n,\,\, r+1), \qquad\text{for $0 \leq r
    \leq j$},
$$
the slope of the  line segment from $\bee_i^{r}$ to $\bee_i^{r+1}$ is $\frac{\Delta x}{\Delta y} = a_{-(s_i+r+1)}$
for $-1 \leq r \leq j$. That is, the $j+1$ segments of $\mathcal{C}$ have slopes $a_{-s_i},a_{-(s_i+1)},
a_{-(s_i+2)}, \ldots, a_{-(s_i+j)}$, listed in order from left to right. Since $\shift{\ba}{s_i+j} =
(a_{-(s_i+j)},a_{-(s_i+j-1)},\ldots,a_{-(s_i+j-m+1)})$, this identifies $\mathcal{C}$ as the first $j+1$ segments
of the boundary curve $\bdry{\shift{\ba}{s_i+j}}$ rotated $180^{\circ}$ and translated.

Shift $\mathcal{C}$ to the left one unit to obtain a new curve $\mathcal{C}'$ that terminates at $(x-1,j+1)$.
(See Figure~\ref{fig:bijproof}B.) Since, by definition, $\pth{P}$ remains weakly below
$\bdry{\shift{\ba}{s_i+j}}$, so too does the subpath $\pth{P}_1$. Since $\pth{P}_1'$ and $\mathcal{C}'$ are
obtained by rotating $\pth{P}_1$ and $\bdry{\shift{\ba}{s_i+j}}$, respectively, it follows that $\pth{P}_1'$
must remain weakly above $\mathcal{C}'$. But $\bee_i^0,\ldots,\bee_i^{j-1}$ lie on $\mathcal{C}$, so they lie
strictly below $\mathcal{C}'$, and therefore $\pth{P}_1'$ avoids these points.  The same is clearly true of
$\pth{P}'$, and this establishes that $\psi_j$ is well defined.

We claim $\psi_j : \setbp_i^j \longrightarrow \setap{j}$ is a bijection. Observe that this establishes
Lemma~\ref{thm:countbadpaths}, since the sets $\setbp_i(\ba,t) = \bigcup_{j=0}^{l-1} \setbp_i^j$ and $\setau
= \bigcup_{i=0}^{l-1} \setau^j$ are then equinumerous and the cardinality of $\setau$ is clearly ${k + l
\choose l-1}$.

\begin{figure}
\begin{center}
    \includegraphics[width=.9\textwidth]{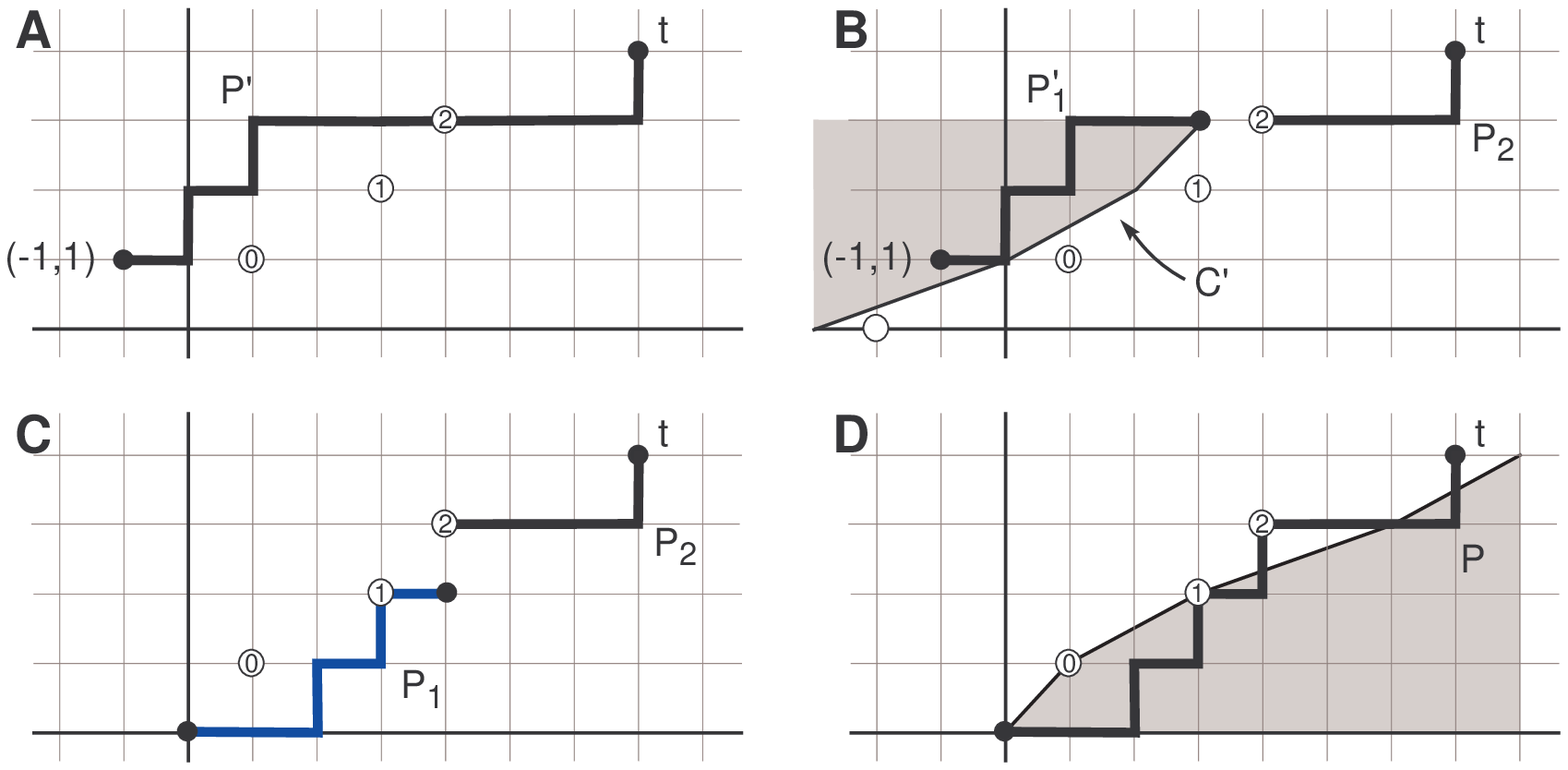}
    \caption{Construction of $\phi_j : \setap{j} \longrightarrow \setbp_i^j$, with
    $\ba=(1,2,3,2)$, $i=4$, $j=2$, and $s_i = 2$. The point $b_i^j$ is indicated with a circled $j$,
    and $b_i^{-1}$ with an open circle.
    The boundary curve in panel D is $\bdry{\shift{\ba}{s_i+j}}=\bdry{\shift{\ba}{4}}$.} \label{fig:inverse}
\end{center}
\end{figure}

To prove $\psi_j$ is bijective, we construct the inverse map $\phi_j : \setap{j} \longrightarrow \setbp_i^j$.
(See Figure~\ref{fig:inverse} for an illustration of the construction.)
Suppose $\pth{P}' \in \setap{j}$.  The first point $\pth{P}'$ intersects amongst $\bee_i^0,\ldots,\bee_i^j$ is
$\bee_i^j$ and it is clear that the step landing at $\bee_i^j$ is horizontal. Remove this step to split $\pth{P}'$
into two paths:  Call the left part $\pth{P}_1'$ and the right part $\pth{P}_2$. Let $\mathcal{C}'$ be the
piecewise linear curve obtained by joining the points $\bee_i^{-1},\bee_i^0,\ldots,\bee_i^j$ (where $\bee_i^{-1}$
is given by~\eqref{eq:bn}) and shifting the result one unit to the left. Then, as above, the segments of
$\mathcal{C}'$ have slopes $a_{-(s_i+j)}, a_{-(s_i+j-1)}, \ldots, a_{-(s_i+1)}, a_{-s_i}$, so that $\mathcal{C}'$
is simply the first $j+1$ segments of $\bdry{\shift{\ba}{s_i+j}}$ rotated $180^\circ$ and translated. Since
$\pth{P}_1'$ lies weakly above $\mathcal{C}'$, the curve $\pth{P}_1$ obtained by rotating $\pth{P}_1'$ by
$180^\circ$ and translating its origin to $(0,0)$ must lie weakly under $\bdry{\shift{\ba}{s_i+j}}$. Attach
$\pth{P}_1$ to $\pth{P}_2'$ by a vertical step to form a new path $\pth{P}$. Then $L=(\pth{P}, (\ba,s_i+j))$ is a
bad LPBP in which $\pth{P}$ terminates at $t$ and has its first bad step landing at $\bee_i^j$. Set
$\phi_j(\pth{P}')=L$, so that clearly $\phi_j \circ \psi_j(L) = L$ and $\psi_j \circ \phi_j(\pth{P}') = \pth{P}'$,
as required. \qed
\newline

Theorem~\ref{thm:gen} now follows immediately from Lemma~\ref{thm:countbadpaths} and the following result:

\begin{lemma}
    \label{thm:completeness}
     Suppose $0 \leq k \leq n$, $0 \leq l \leq m$, and the point
    $(k+1,l)$ is weakly right of $\bdry{\shift{\ba}{j}}$ for all $j$.
    Then each of the sets $\Bee_0,\ldots,\Bee_{n-1}$ is complete with respect to $t=(k,l)$.
\end{lemma}

\noindent \textbf{Proof.} The set $\Bee_i$ is not complete with respect to $t$ if and only if it contains some
point of the form $(k+\delta,l)$ with $\delta \geq 1$.  But $\Bee_i$ consists of those points that lie immediately
above $\bdry{\shift{\ba}{j}}$ for some $j$. Thus $(k+\delta,l) \in \Bee_i$ for some $\delta \geq 1$ precisely when
$(k+1,l)$ is strictly left of some boundary $\bdry{\shift{\ba}{j}}$. The result follows. \qed

\section{A Cycle Lemma Proof of Theorem~\ref{thm:gen}}
\label{sec:cyclelemma}

We now sketch an alternative proof of Theorem~\ref{thm:gen} using the cycle lemma~\cite{dmo}.  The
formulation most applicable here is the following:\newline

\noindent\textbf{Cycle Lemma.} \emph{Let $\mathbf{i}=(i_0,\ldots,i_m)$ be a sequence with integral entries
$i_j \leq 1$ having positive sum $k = i_0 + \cdots + i_m$. Then there are exactly $k$ cyclic shifts of
$\mathbf{i}$ with all partial sums positive.} \qed \newline

\newcommand{\RR}{\mathsf{R}}
\newcommand{\UU}{\mathsf{U}}
\newcommand{\bt}{\mathbf{u}}
\newcommand{\good}{\mathscr{G}^*(\ba,t)}
\newcommand{\words}{\mathscr{W}}
\newcommand{\cycmap}{\Omega}

The following result is the key to our alternative proof of Theorem~\ref{thm:gen}.

\begin{lemma} \label{lem:cycleproof}
Let $\ba$ be a weak $m$-part composition of $n$, and let $t=(n,l)$, where $0 \leq l < m$.  Let $\good$ be the
set of good LPBPs of the form $(\pth{P},(\ba,j))$ where $\pth{P}$ is a path from $(0,0)$ to $t$ that
terminates with a right step. Then
$$
        |\good| =  \binom{n+l-1}{l}(m-l).
$$
\end{lemma}

\noindent \textbf{Proof:}  Let $\words$ be the set of words of length $n+l$ on the alphabet $\{\RR,\UU\}$
that contain $l$ $\UU$'s and $n$ $\RR$'s and end with an $\RR$. We give a bijection $\cycmap \,:\, \words
\times [m-l] \longrightarrow \good$, where  $[m-l]=\{1,\ldots,m-l\}$.  The construction is illustrated in
Example~\ref{exmp:cycleproof}, below.

Let $(w,k) \in \words \times [m-l]$.  Factor $w$ into $m$ blocks $w=w_0\cdots w_{m-1}$ as follows:  Suppose
$\ba=(a_0,\ldots,a_{m-1})$, and parse $w$ from left to right letting $w_0,\ldots,w_{m-1}$ in turn be maximal
contiguous substrings such that
\begin{itemize}
\item   $w_i$ is empty if $a_i=0$,
\item   $w_i$ contains $a_i$ $\RR$'s and ends with an $\RR$ if $a_i > 0$.
\end{itemize}
Observe that this decomposition of $w$ is unique.

Consider the integer sequence $\bt=(1,-u_0,1,-u_1,\ldots,1,-u_{m-1})$, where $u_i$ is the number of $\UU$'s
in $w_i$.  The entries of $\bt$ sum to $m-(u_0+\ldots+u_{m-1})=m-l > 0$, so the cycle lemma implies there are
exactly $m-l$ cyclic shifts of $\bt$ whose partial sums are all positive. Clearly such shifts must be of the
form $\shift{\bt}{-2s}$, where $0 \leq s < m$. (See~\eqref{eq:shift} for the definition of
$\shift{\bt}{-2s}$.) Suppose the good shifts are $\shift{\bt}{-2s_1},\ldots,\shift{\bt}{-2s_{m-l}}$, where
$s_1 < \cdots < s_{m-l}$. Set $j=s_{k}$ and form the word $w'=w_{j} w_{j+1} \cdots w_{j+m-1}$, where the
indices are to be interpreted modulo $m$. From $w'$, construct a lattice path $\pth{P}$ originating at
$(0,0)$  by treating $\RR$ and $\UU$ as right and up steps, respectively.

Set $\cycmap(w,k) :=  (\pth{P}, (\ba,-j))$. Observe that indeed $\cycmap(w,k) \in
\good$, since $\pth{P}$ clearly terminates at $(n,l)$ with a right step and
\begin{align*}
            &\text{ $\shift{\bt}{-2j}$ has all partial sums positive} \\
    \Longleftrightarrow\quad &\text{ $u_j+ \cdots + u_{j+d} < d+1$, for $0 \leq d < m$} \\
    \Longleftrightarrow\quad &\text{$\pth{P}$ has at least $a_j+\cdots+a_{j+d}$ right steps before its $(d+1)$-st
    up step,
    for $0 \leq d < m$} \\
    \Longleftrightarrow\quad &\text{$\pth{P}$ is dominated by $\shift{\ba}{-j}=(a_{j},a_{j+1},\ldots,a_{j+m-1})$.}
\end{align*}
Moreover, this construction of $(\pth{P},(\ba,-j))$ from $(w,k)$ can be reversed, as follows: (1) recover
$w'$ from $\pth{P}$, (2) parse $w'$ as above, but relative to the composition $\shift{\ba}{-j}$, to obtain
$w_j,w_{j+1},\ldots,w_{j+d}$ and hence $w$, (3) retrieve $\bt$ from  $w_0,\ldots,w_{m-1}$, and (4) deduce $k$
by applying the cycle lemma to $\bt$.

Thus $\cycmap \,:\, \words \times [m-l] \longrightarrow \good$ is bijective, and since $|\words|=
\binom{n+l-1}{l}$, the result follows. \qed

\begin{figure}
\begin{center}
    \includegraphics[width=.35\textwidth]{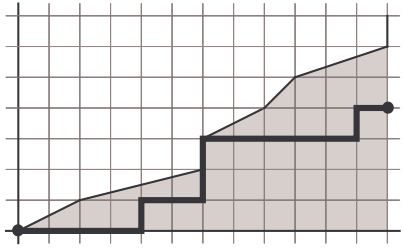}
    \caption{The path constructed in Example~\ref{exmp:cycleproof}\label{fig:cycle}.}
\end{center}
\end{figure}

\begin{example} \label{exmp:cycleproof}
Let $n=12$, $m=7$, $l=4$, $\ba=(1,3,0,2,4,0,2)$ and take
$$
    (w,k)=(\RR\RR\RR\UU\RR\RR\RR\RR\RR\UU\RR\RR\UU\UU\RR\RR, 3).
$$
Then we have
$$
    w_0=\RR, \quad
    w_1=\RR\RR\UU\RR, \quad
    w_2=\epsilon, \quad
    w_3=\RR\RR, \quad
    w_4=\RR\RR\UU\RR\RR, \quad
    w_5=\epsilon, \quad
    w_6=\UU\UU\RR\RR,
$$
where $\epsilon$ denotes the empty string.  This gives $\bt=(1,0,1,-1,1,0,1,0,1,-1, 1,0,1,-2)$, and the
$m-l=3$ cyclic shifts of $\bt$ with all partial sums positive are seen to be
\begin{align*}
    \shift{\bt}{0} &=(1,0,1,-1,1,0,1,0,1,-1, 1,0,1,-2) \\
    \shift{\bt}{-4} &= (1,0,1,0,1,-1,1,0,1,-2,1,0,1,-1) \\
    \shift{\bt}{-6} &= (1,0,1,-1,1,0,1,-2,1,0,1,-2,1,0).
\end{align*}
Thus $s_1 = 0$, $s_2=2$, $s_3=3$, so that $j=s_3 = 3$ and
$$
    w' = w_3 w_4 w_5 w_6 w_0 w_1 w_2 = \RR\RR\RR\RR\UU\RR\RR\UU\UU\RR\RR\RR\RR\RR\UU\RR.
$$
Figure~\ref{fig:cycle} shows the path $\pth{P}$ corresponding to $w'$ and the dominating boundary
$\bdry{\shift{\ba}{-j}}=\bdry{\shift{\ba}{-3}}$. \qed
\end{example}

It is now easy to establish Theorem~\ref{thm:gen} in the special case where the terminal point is $t=(n,l)$
for some $0 \leq l \leq m$. In particular, we clearly have $|\setgg{(n,l)}{\ba}| = \sum_{i=0}^{l}
|\mathscr{G}^*(\ba,(n,i))|$, so Lemma~\ref{lem:cycleproof} gives
\begin{align}
\label{eq:endpoint}
    |\setgg{(n,l)}{\ba}| =
    \sum_{i=0}^{l} (m-i) \binom{n+i-1}{i}
                = m\binom{n+l}{l} - n\binom{n+l}{l-1},
\end{align}
in agreement with~\eqref{eq:breakdown}.  Moreover, we have the usual lattice path recursion
$$
    |\setgg{(k+1,l)}{\ba}| = |\setgg{(k,l)}{\ba}| + |\setgg{(k+1,l-1)}{\ba}|
$$
provided $(k+1,l)$ is weakly right of every shift of $\bdry{\ba}$. Using~\eqref{eq:endpoint} as an initial
condition and iterating the above recursion allows us to determine $|\setgg{(k,l)}{\ba}|$ for any terminal point
$(k,l)$ satisfying this same condition. Thus we have an inductive proof of~\eqref{eq:breakdown}. The details are
more tedious than illuminating and are omitted here.

\section{A Refinement:  Counting Paths with a Specified Number of Corners}
\label{sec:refinements}

\newcommand{\setcc}[1]{\mathscr{C}^{#1}}
\newcommand{\setupc}[1]{\mathscr{C}_{#1}}

An \emph{up-right corner} in a lattice path is a point at which an up step terminates and is immediately
followed immediately by a right step. Observe that for any $c$ points $(X_j,Y_j)$ satisfying
$$
    (0,1) \wlwb (X_1,Y_1) \lb (X_2,Y_2) \lb \cdots \lb (X_c,Y_c) \wlwb (k-1,l),
$$
there is a unique lattice path from $(0,0)$ to $(k,l)$ having up-right corners at exactly these points. Since
the the $X_j$'s and $Y_j$'s can be chosen in $\binom{k}{c}$ and $\binom{l}{c}$ ways, respectively, it follows
that there are ${k \choose c}{l \choose c}$ lattice paths from $(0,0)$ to $(k,l)$  with exactly $c$ up-right
corners.

Let $\setcc{c}$ be the set of all LPBPs whose  paths have exactly $c$ up-right corners.   The following
theorem is a generalization of~\cite[Theorem 3.4.2]{kratt:1}, and our proof is inspired by that
of~\cite[Theorem 5]{goulserr:1}.

\begin{theorem}
    \label{thm:corners}
    Let $\ba$ be any weak $m$-part composition of $n$ and let $t=(k,l)$ be a point
    dominated by all cyclic shifts of $\ba$.  Then
    \begin{equation*}
        |\setgg{t}{\ba} \cap \setcc{c}| = m {k \choose c}{l \choose c} - n  {k-1 \choose c-1} {l+1 \choose c+1}.
    \end{equation*}
    \qed
\end{theorem}

Clearly this is a refinement of Theorem~\ref{thm:gen}, and unsurprisingly our proof relies on a corresponding
refinement of Lemma~\ref{thm:countbadpaths}.  Note that the hypothesis regarding the terminal point $t=(k,l)$
is slightly stronger than that of Theorem~\ref{thm:gen}. That is, we require $(k,l)$, rather than $(k+1,l)$,
to be dominated by all $\shift{\ba}{j}$.

Indeed, our refinement of Lemma~\ref{thm:countbadpaths} requires a slightly stronger notion than
completeness.  With the sets $\Bee_i$ defined as in~\eqref{eq:setbi}, we say $\Bee_i$ is \emph{strongly
complete} with respect to $t=(k,l)$ if the point $b_i^{l-1}$ is \emph{strictly} to the left of $t$.

Theorem~\ref{thm:corners} follows immediately from the following two results. We assume the notation of
Section~\ref{sec:proof} throughout.

\begin{lemma}
    \label{thm:strongcompleteness}
     Suppose $0 \leq k \leq n$, $0 \leq l \leq m$, and the point
    $(k,l)$ is dominated by $\shift{\ba}{j}$ for all $j$.
    Then each of the sets $\Bee_0,\ldots,\Bee_{n-1}$ is strongly complete with respect to $t=(k,l)$.
\end{lemma}

\noindent\textbf{Proof.}  This is an obvious modification of Lemma~\ref{thm:completeness}. \qed

\begin{lemma}\label{thm:refine}
    If $\Bee_i$ strongly complete with respect to $t=(k,l)$, then
    \begin{equation*}
        |\setbp_i(\ba,t) \cap \setcc{c}| =  {k-1 \choose c-1} {l+1 \choose c+1}.
    \end{equation*}
\end{lemma}

\newcommand{\bx}{\mathbf{X}}
\newcommand{\by}{\mathbf{Y}}

\noindent \textbf{Proof.}  We prove the lemma by giving a bijection between $\setbp_i(\ba,t) \cap \setcc{c}$ and
pairs of sequences $(\bx,\by) \in \mathbb{Z}^{c} \times \mathbb{Z}^{c+1}$ satisfying
\begin{equation*}
\label{eq:seq}
    0 \leq X_1 < \cdots < X_c = k-1 \qquad\text{and}\qquad 1 \leq Y_1 < \ldots < Y_{c+1} \leq l+1.
\end{equation*}
Fix such a pair $(\bx,\by)$. Since $\Bee_i$ is strongly complete with respect to $(k,l)$, we have
$$x(\bee_i^{c-1}) \leq k-1 = X_c,$$ where $x(p)$ denotes the $x$-coordinate of the point $p$.  Let $r \leq c$ be
the smallest index for which $x(\bee_i^{r-1}) \leq X_r$, and set $j=Y_r-1$ so that  $\bee_i^j =
(x(\bee_i^j),Y_r)$. Since $Y_r \geq r$, we have $j \geq r-1$, so the minimality of $r$ implies either $r=1$ or
$X_{r-1} < x(b_i^{r-2}) \leq x(b_i^j)$. Thus we have a chain of points
\begin{align*}
    (X_1, Y_1) \lb \cdots \lb (X_{r-1}, Y_{r-1}) \lb\; &\bee_i^j
        \wlwb
    (x(\bee_i^j), Y_{r+1}-1)
        \lb \\
    &(X_{r}+1, Y_{r+2} - 1) \lb \cdots \lb (X_{c-1}+1, Y_{c+1} - 1).
\end{align*}
It is easy to verify that there is a unique path $\pth{P}$ from $(-1,1)$ to $(k,l)$ passing through all these
points such that:
\begin{itemize}
\item   $\pth{P}$ has $r-1$ right-up corners at $(X_1, Y_1) \ldots (X_{r-1}, Y_{r-1})$,
        and no further right-up corners strictly left of $\bee_i^j$,
\item   the steps of $\pth{P}$ terminating at $\bee_i^j$ and originating at $(x(\bee_i^j), Y_{r+1}-1)$ are
horizontal,
\item   $\pth{P}$ has $c-r$ up-right corners at $(X_{r}+1, Y_{r+2} - 1), \ldots, (X_{c-1}+1, Y_{c+1} - 1)$,
        and no further up-right corners strictly right of $\bee_i^j$.
\end{itemize}

By construction,  $\pth{P}$ avoids $\bee_i^0,\ldots,\bee_i^{j-1}$ but meets $\bee_i^j$. So we can apply the
bijection $\phi_j$ (see the proof of Lemma~\ref{thm:countbadpaths}) to get an LPBP $\phi_j(\pth{P}) = (\pth{P'},
(\ba,s_i+j)) \in \setbp_i(\ba,t)$. Observe that the $r-1$ right-up corners of $\pth{P}$ to the left of $\bee_i^j$
become up-right corners of $\pth{P}'$ through rotation, while the $c-r$ up-right corners of $\pth{P}$ to the right
of $\bee_i^j$ are preserved in $\pth{P}'$.

We now check for corners at $\bee_i^j$ and $(x(\bee_i^j), Y_{r+1}-1)$.  There are two cases to consider. If
$Y_{r+1} - 1 > Y_{r}$, then $\pth{P}'$ does not have an up-right corner at $\bee_i^j$ but does at
$(x(\bee_i^j), Y_{r+1}-1)$.  Otherwise $Y_{r+1} -1 = Y_r$, in which case $\bee_i^j = (x(\bee_i^j),
Y_{r+1}-1)$ and $\pth{P}'$ has an up-right corner at this point.

In either case, $\pth{P}'$ has exactly $(r-1)+(c-r)+1=c$ up-right corners in total. That is, $\phi_j(\pth{P})
\in \setbp_i(\ba,t) \cap \setcc{c}$. Since $\phi_j$ is bijective, so too is the correspondence $(\bx,\by)
\mapsto \phi_j(\pth{P})$ described here. This completes the proof.\qed
\\

In analogy with up-right corners, we say a \emph{right-up corner} is formed when a right step is followed
immediately by an up step.  It is convenient to treat an initial up step as a \emph{virtual} right-up corner.
Then, letting $\setupc{c}$ be the set of all LPBPs whose paths have exactly $c$ right-up corners (real or
virtual), we have:

\begin{theorem}
    \label{thm:corners2}
    Let $\ba$ be any weak $m$-part composition of $n$ and let $t=(k,l)$ be a point
    dominated by all cyclic shifts of $\ba$.  Then
    \begin{equation*}
        |\setgg{t}{\ba} \cap \setupc{c}| = m {k+1 \choose c}{l-1 \choose c
        -1} - n  {k \choose c-1} {l \choose c}.
    \end{equation*}
    \qed
\end{theorem}

Consider the case $(k,l)=(n,m)$ in Theorems~\ref{thm:corners} and~\ref{thm:corners2}.  Notice that the first and
last corners of any good path are right-up corners.  Since right-up corners and up-right corners must alternate,
the number of good paths with $c$ right-up corners is equal to the number of good paths with $c-1$ up-right
corners. Indeed, Theorems \ref{thm:corners} and \ref{thm:corners2} show this common number to be ${n \choose c-1}
{m \choose c}$.

\section{Counting Paths Dominated by Periodic Boundaries}
\label{sec:periodic}

\newcommand{\abn}[3]{(#1,#2)^{#3}}
\newcommand{\ab}[2]{\partial_{#1,#2}}
\newcommand{\stair}[2]{\pth{S}_{#1,#2}}
\newcommand{\Dom}[2]{D(#1,#2)}

Let $\ba=(a_0,\ldots,a_{m-1})$ be a weak composition of $n$, and let $d$ be the least positive integer such that
$\shift{\ba}{d} = \ba$. Clearly $d$ divides $m$. In the case that $d < m$ we say $\ba$ is \emph{periodic} with
\emph{period $d$}.  For example $(3,1,2,3,1,2)$ has period 3.

If $\ba$ has period $d$ then $\shift{\ba}{i} = \shift{\ba}{i\!\mod d}$.  Thus
\begin{equation}
\label{eq:periodic}
    |\setgg{t}{\ba}| =
    \frac{m}{d}\big(\Dom{\ba}{t} + \Dom{\shift{\ba}{1}}{t} + \cdots + \Dom{\shift{\ba}{d-1}}{t}\big),
\end{equation}
where $\Dom{\ba}{t}$ denotes the number of paths from $(0,0)$ to point $t$ dominated by $\ba$.  The left-hand side
of this equality can generally be evaluated by Theorem~\ref{thm:gen}, and in certain special cases this allows us
to deduce $\Dom{\ba}{t}$.

The case $d=1$ is particularly straightforward. Here we have $\ba=(a)^m=(a,a,\ldots,a)$,
and~\eqref{eq:periodic} gives $\Dom{\ba}{t}=\frac{1}{m}|\setgg{t}{\ba}|$.  This is precisely our earlier
proof of  Corollary~\ref{cor:kratt}.

\begin{figure}
    \includegraphics[width=.35\textwidth]{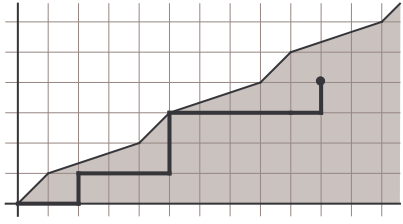}
    \caption{The boundary $\ab{1}{3}$ and a path it dominates.}\label{fig:boundaries}
\end{figure}

The case $d=2$ involves compositions of the form $\ba=(a,b,a,b,\ldots,a,b)=(a,b)^m$, where $a,b$ are distinct
nonnegative integers. In the following discussion it will be convenient to write $\ab{a}{b}$ for the ``infinite''
boundary curve $\bdry{(a,b,a,b,\ldots)}$. See Figure~\ref{fig:boundaries} for an illustration of $\ab{1}{3}$ and
one path that it dominates.

Theorem~\ref{thm:periodic}, below, gives explicit formulae for the number of paths under $\ab{a}{b}$ or
$\ab{b}{a}$ to certain special endpoints. This result is also implicit in Tamm~\cite[Propositions 2,3]{tamm},
where it appears in generating series form. The proof given there follows a probabilistic argument (originally due
to Gessel) reliant on Lagrange inversion, whereas our derivation is purely combinatorial.

\newcommand{\pP}[1]{p_{#1}}
\newcommand{\pQ}[1]{q_{#1}}
\newcommand{\baab}{{\mathbf{a}}}
\newcommand{\baba}{{\mathbf{b}}}
\newcommand{\psetPa}[1]{\mathscr{P}^{a,b}_{#1}}
\newcommand{\psetPb}[1]{\mathscr{P}^{b,a}_{#1}}
\newcommand{\psetQa}[1]{\mathscr{Q}^{a,b}_{#1}}
\newcommand{\psetQb}[1]{\mathscr{Q}^{b,a}_{#1}}
\newcommand{\Mab}[1]{M_{#1}}
\newcommand{\Nab}[1]{N_{#1}}
\newcommand{\psetP}[3]{\mathscr{P}^{#1,#2}_{#3}}
\newcommand{\psetQ}[3]{\mathscr{Q}^{#1,#2}_{#3}}

\begin{theorem}
\label{thm:periodic} Fix integers $a,b$ with $0 \leq a < b$ and set $c=a+b$. For $n \geq 0$ let
\begin{align*}
    \pP{n} = (cn+b-a-1,2n) \qquad\text{and}\qquad
    \pQ{n} = (cn+b-1,2n+1).
\end{align*}
Let $\psetPa{n}$ and $\psetQa{n}$, respectively, be the sets of lattice paths from the origin to $\pP{n}$ and
$\pQ{n}$ that lie weakly under $\ab{a}{b}$.  Define sets $\psetPb{n}$ and $\psetQb{n}$ similarly, but for paths
weakly under $\ab{b}{a}$. Then
\begin{xalignat}{2}
\label{eq:per1}
    |\psetQa{n}| &= \Mab{n},
    &
    |\psetQb{n}| &= 0, \\
\label{eq:per2}
    |\psetPa{n}| &= \Nab{n}
                    + \frac{1}{2} \sum_{i=0}^{n-1} \Mab{i}\Mab{n-1-i},
    &
    |\psetPb{n}| &= \Nab{n}
                    - \frac{1}{2} \sum_{i=0}^{n-1} \Mab{i}\Mab{n-1-i}.
\end{xalignat}
where
\begin{align*}
    \Mab{n} = \frac{b-a}{cn+b}\binom{ (c+2)n + b}{ 2n +1 } \qquad\text{and}\qquad
    \Nab{n} = \frac{b-a}{cn+b-a} \binom{(c+2)n + b-a-1}{ 2n }.
\end{align*}
\end{theorem}

\begin{figure}
    \includegraphics[width=.900\textwidth]{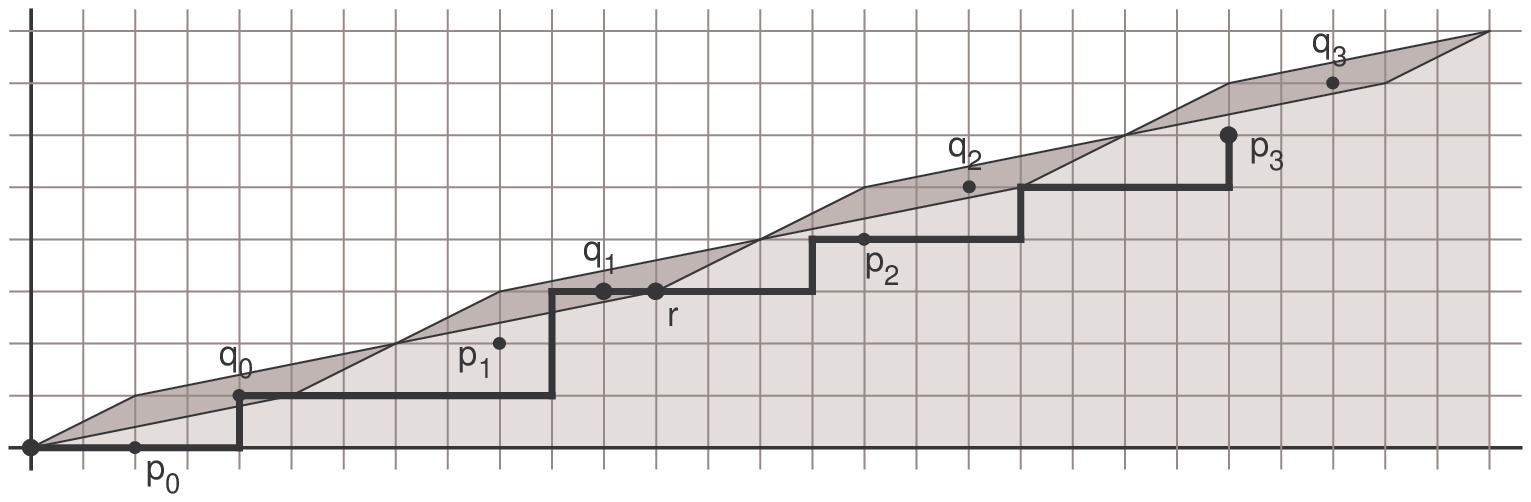}
    \caption{The points $\pP{i},\pQ{i}$ in the case $a=2,b=5$ and a
path $\pth{P} \in \psetP{2}{5}{3}$.}\label{fig:pnqn2}
\end{figure}

\noindent \textbf{Proof.} A glance at Figure~\ref{fig:pnqn2} will make the proof more clear.  It illustrates
several points $\pP{i},\pQ{i}$ in the case $a=2,b=5$, along with the boundaries $\ab{2}{5}$, $\ab{5}{2}$ and a
path $\pth{P} \in \psetP{2}{5}{3}$.

Let $\baab=(a,b)^{n+1}$ and $\baba=(b,a)^{n+1}$, so a path lies under $\ab{a}{b}$ (respectively, $\ab{b}{a}$) if
and only if it is dominated by $\baab$ (respectively, $\baba$).

Clearly $|\psetQb{n}|=0$ since $\pQ{n}$ is not dominated by $\ab{b}{a}$. Moreover, since $\shift{\baab}{2j}=\baab$
and $\shift{\baab}{2j+1}=\baba$ for all $j$, we have
\begin{equation*}
    |\setgg{\pQ{n}}{\ba}| = (n+1)\big(|\psetQa{n}|+|\psetQb{n}|\big) = (n+1) |\psetQa{n}|.
\end{equation*}
The point $r:=(cn+b,2n+1)$ one unit right of $\pQ{n}$ is dominated by both $\baab$ and $\baba$.
Theorem~\ref{thm:gen} may therefore be applied, and it gives $|\setgg{\pQ{n}}{\ba}|=(n+1)\Mab{n}$. This
establishes~\eqref{eq:per1}. A similar analysis yields
\begin{align}
\label{eq:pth1}
    |\psetPa{n}|+|\psetPb{n}| &= \frac{1}{n+1}|\setgg{\pP{n}}{\ba}|
                                = 2\Nab{n}.
\end{align}

Observe that $\psetPb{n} \subset \psetPa{n}$. In fact, a path is dominated by $\baba$ if and only if it is
dominated by $\baab$ and misses each of the points $\pQ{0},\ldots,\pQ{n}$.  Now consider a path $\pth{P} \in
\psetPa{n} \setminus \psetPb{n}$, and let $i$ be the largest index so that $\pth{P}$ meets $\pQ{i}$. Then
$\pth{P}$ exits $\pQ{i}$ with a right-step to $r$, and removal of this step splits $\pth{P}$ into two paths,
$\pth{P}'$ and $\pth{P}''$, with $\pth{P}' \in \psetQa{i}$ and $\pth{P}''$ dominated by $\ab{a}{b}$. (See
Figure~\ref{fig:pnqn2}.) But $\pP{n}-r=\pQ{n-i-1}$, so we effectively have $\pth{P}'' \in \psetQa{n-1-i}$. Hence
$$
    |\psetPa{n}| = |\psetPb{n}| + |\psetPa{n} \setminus \psetPb{n}|
            = |\psetPb{n}| + \sum_{i=0}^{n-1} |\psetQa{i}| |\psetQa{n-i-1}|.
$$
Formulae~\eqref{eq:per2} now follow from~\eqref{eq:per1} and~\eqref{eq:pth1}. \qed \\

\newcommand{\overc}{\overline{c}}
\newcommand{\underc}{\underline{c}}

When $a=\frac{c-1}{2}$, $b=\frac{c+1}{2}$ for an odd positive integer $c$, observe that a path is dominated
by $\ab{b}{a}$ if and only if it lies weakly under the line $cx=2y$. So Theorem~\ref{thm:periodic} can be
applied to give the following enumeration of paths under a line of half-integer slope. An equivalent result
also appears as~\cite[Theorem 1]{tamm}.

\begin{corollary}
\label{cor:periodic1} Let $c$ be an odd positive integer. The number of lattice paths from $(0,0)$ to $(cn,2n)$
that lie weakly below the line $cx=2y$ is given by
$$
    \frac{1}{cn+1} \binom{(c+2)n}{2n}
        - \frac{1}{2} \sum_{i=0}^{n-1} \Mab{i} \Mab{n-1-i},
        \qquad\text{where}\quad
        \Mab{i} = \frac{1}{2i+1} \binom{(c+2)i+\frac{c+1}{2}}{2i}.
$$
\qed
\end{corollary}

\begin{figure}
    \includegraphics[width=.750\textwidth]{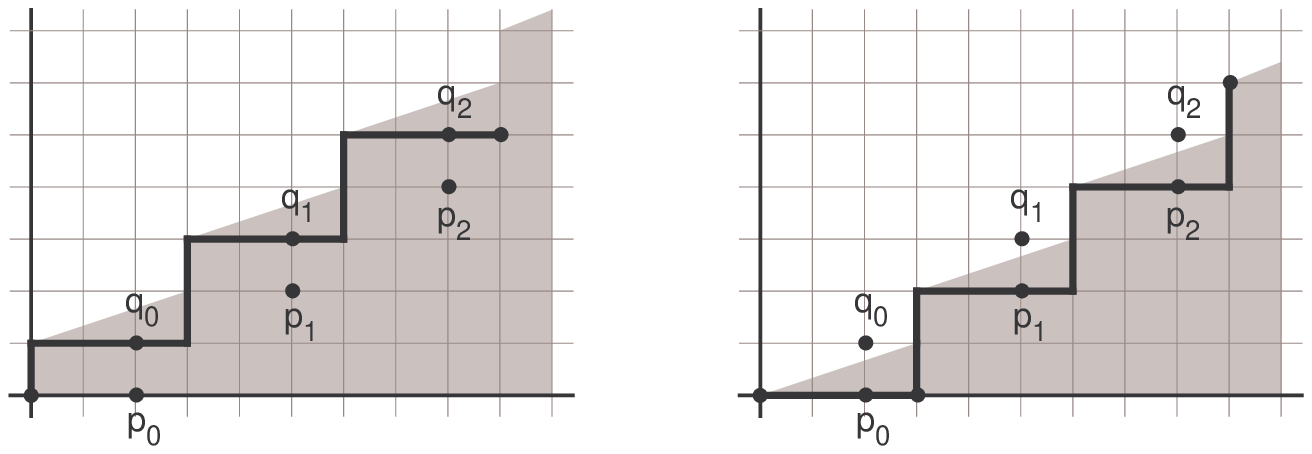}
    \caption{The boundaries $\ab{0}{3}$ and $\ab{3}{0}$ and corresponding staircases
    $\UU (\RR^3 \UU^2)^2 \RR^3$ and $(\RR^3 \UU^2)^3$.}\label{fig:staircase}
\end{figure}

Another special case of Theorem~\ref{thm:periodic} worth mentioning is that when $a=0$, where we count paths from
$(0,0)$ to $\pP{n}=(b(n+1)-1,2n)$ or $\pQ{n}=(b(n+1)-1,2n+1)$ dominated by $\ab{b}{0}$ or $\ab{0}{b}$. Observe
that a path to either point is bounded by $\ab{0}{b}$ if and only if it lies weakly under the ``staircase'' $\UU
(\RR^b \UU^2)^n \RR^b$. (See Figure~\ref{fig:staircase}.)
Similarly, the paths dominated by $\ab{b}{0}$ are precisely those that lie weakly under $(\RR^b \UU^2)^n$. For $n
\geq 1$, such paths begin with $b$ right steps, and removing these puts $\psetP{b}{0}{n}$ in correspondence with
the set of paths from $(0,0)$ to $(bn-1,2n)$ that lie weakly beneath $(\UU^2 \RR^b)^n$.

When $a=0,b=2$, the various quantities in Theorem~\ref{thm:periodic} can be compactly expressed in terms of the
Catalan numbers. In particular, we obtain simple formulae for the number of paths from $(0,0)$ to any point on the
boundary $\UU(\RR^2\UU^2)^n\RR$.  Note that the usual recursions for lattice paths then give similar expressions
for paths to any point near the boundary.

\begin{corollary}
\label{cor:periodic2}  Let $C_n = \frac{1}{n+1}\binom{2n}{n}$ be the $n$-th Catalan number. There are
\begin{itemize}
\item   $2C_{2n+1}$ paths from $(0,0)$ to $(2n+1,2n+1)$, and

\item   $2^{2n+1}C_n-C_{2n+1}$ paths from $(0,0)$ to $(2n,2n)$ or $(2n,2n \pm 1)$
\end{itemize}
that lie weakly under $\UU(\RR^2\UU^2)^n\RR$. Moreover, there are  $2^{2n+1}C_n-C_{2n+1}$  paths from $(0,0)$ to
$(2n-1,2n)$ lying weakly under $(\UU^2\RR^2)^n$.
\end{corollary}

\noindent \textbf{Proof.}  The desired result rests upon the convolution identity
\begin{equation}
\label{eq:conv}
    \sum_{i=0}^{n-1} C_{2i+1} C_{2n-(2i+1)} = C_{2n+1}-2^{2n}C_n.
\end{equation}
This is easily seen to be equivalent to the functional equation
\begin{equation}
\label{eq:cat}
    D(x)^2=\frac{1}{x}D(x)-C(4x^2),
\end{equation}
where $C(x) = \sum_n C_n x^n$ is the Catalan generating series and $D(x)=\frac{1}{2}(C(x)-C(-x))$ is its odd part.
To establish~\eqref{eq:cat}, expand $D(x)^2 = \frac{1}{4}(C(x)-C(-x))^2$ and substitute
\begin{align*}
    C(x)^2+C(-x)^2 &= \tfrac{2}{x}D(x) \\
    C(x)C(-x) &= 2C(4x^2)-\tfrac{1}{x}D(x),
\end{align*}
which themselves are readily derived from the well-known identities  $C(x)=1+xC(x)^2$ and $C(x) =
\frac{1}{2x}(1-\sqrt{1-4x})$, respectively.

Apply Theorem~\ref{thm:periodic} with $a=0,b=2$, noting that $\Mab{n}=2C_{2n+1}$, $\Nab{n}=C_{2n+1}$ and
using~\eqref{eq:conv} to simplify the results. This gives $|\psetQ{0}{b}{n}| = 2C_{2n+1}$ paths to $(2n+1,2n+1)$
and  $|\psetP{0}{b}{n}| = 3C_{2n+1}-2^{2n+1}C_n$ paths to $(2n+1,2n)$ under $\UU(\RR^2\UU^2)^n\RR$. Since paths to
$(2n+1,2n+1)$ pass through either $(2n+1,2n)$ or $(2n,2n)$, there are
$|\psetQ{0}{b}{n}|-|\psetP{0}{b}{n}|=2^{2n+1}-C_{2n+1}$ paths to $(2n,2n)$ under $\UU(\RR^2\UU^2)^n\RR$. Clearly
there are this same number of paths to $(2n,2n\pm 1)$.

Finally, paths to $(2n-1,2n)$ under $(\UU^2\RR^2)^n$ are in bijection with $\psetP{b}{0}{n}$, and
Theorem~\ref{thm:periodic} yields $|\psetP{b}{0}{n}| = 2^{2n+1}C_n-C_{2n+1}$. Alternatively, we could rotate and
flip to view these as paths from $(0,0)$ to $(2n,2n-1)$ dominated by $\UU(\RR^2\UU^2)^n\RR$. \qed
\\

Let $\ba$ be a composition of period $d$, and consider a terminus $t=(k,l)$ such that the point $(k+1,l)$ is
dominated by all cyclic shifts of $\ba$, but  no shift except $\ba$ itself dominates $t$. Then we clearly have
$\Dom{\shift{\ba}{i}}{t} = 0$ for $i \geq 1$, so applying Theorem~\ref{thm:gen} in tandem with~\eqref{eq:periodic}
gives a closed form expression for $\Dom{\ba}{t}$. Indeed, the key to our proof of Theorem~\ref{thm:periodic} was
to determine $|\psetQa{n}|$ in exactly this way.

As another interesting example we present the following result, also recently discovered independently by other
authors~\cite{chow}. (It appears there in a very slightly modified form. We shall make further comments below.)

\begin{theorem}
\label{thm:periodic3} Let $s,t$ and $n$ be positive integers.  Then there are
$$
\frac{1}{n}\binom{(s+t)n-2}{tn-1}
$$
lattice paths from $(0,0)$ to $(sn-1,tn-1)$ lying weakly beneath $\UU^{t-1}(\RR^s\UU^t)^{n-1}\RR^{s-1}$.
\end{theorem}

\noindent \textbf{Proof.}  Let $\ba = (0^{t-1},s)^n$, so that $\ba$ is a $tn$-part composition of $sn$ with period
$t$. Note that a path from $(0,0)$ to $(sn-1,tn-1)$ lies weakly beneath $\UU^{t-1}(\RR^s\UU^t)^{n-1}\RR^{s-1}$
precisely when it is dominated by $\ba$.  Furthermore, none of $\shift{\ba}{1},\ldots,\shift{\ba}{t-1}$ dominate
$(sn-1,tn-1)$, whereas all of them dominate $(sn,tn-1)$.  The result follows immediately from~\eqref{eq:periodic}
after applying Theorem~\ref{thm:gen} with terminus $(k,l)=(sn-1,tn-1)$. \qed
\\

Setting $s=t=k$ in Theorem~\ref{thm:periodic3} yields the following elegant Catalan result, first appearing
as~\cite[Theorem 8.3]{noy:matroid} with a proof based on the Cycle Lemma. Our need for the terminal point to be
dominated by exactly one cyclic shift of the boundary sheds light on the observation of those authors that the
ostensibly similar problem of counting paths to $(nk,nk)$ dominated by $(\UU^k\RR^k)^n$ is in fact much more
complicated.\footnote{Noy and de Mier~\cite{noy:tennisball} have recently introduced a very elegant approach to
the enumeration of lattice paths from $(0,0)$ to $(sn,tn)$ dominated by $(\UU^t \RR^s)^n$, for arbitrary $s,t$.
They deduce generating series that are products of the fractional power series solutions of a certain functional
equation dependent on $s$ and $t$.}

\begin{corollary}
\label{cor:periodic4} Let $n$ and $k$ be positive integers.  Then there are $kC_{nk-1}$ lattice paths from $(0,0)$
to $(nk-1,nk-1)$ lying weakly beneath $\UU^{k-1}(\RR^k\UU^k)^{n-1}\RR^{k-1}$. \qed
\end{corollary}

We conclude with some comments on recent work by Chapman \etal~\cite{chow}.  They consider lattice paths that
remain strictly below the staircase boundary $\stair{s}{t}$ beginning at $(0,t)$, moving to the right $s$ steps,
then up $t$ steps, to the right $s$ steps, \emph{etc.} That is, $\stair{s}{t}$ is described by $(\RR^s \UU^t)^n$,
but is shifted $t$ units upward to originate at $(0,t)$. Their main results concern the enumeration of two types
of paths avoiding $\stair{s}{t}$, namely those from $(0,0)$ to $(sn + 1, tn)$, and those from $(1,0)$ to $(sn,
tn-1)$. They employ a Cycle Lemma argument similar in structure to our proof of Lemma~\ref{lem:cycleproof} to
obtain compact expressions counting both types of paths, even allowing for the refined enumeration of paths with a
specified number of corners. These same results can be obtained from our methods, as follows.

First observe that a path from $(1,0)$ to $(sn,tn-1)$  avoiding $\stair{s}{t}$ can be shifted left one unit to
give a path from $(0,0)$ to $(sn-1,tn-1)$ lying weakly below $\UU^{t-1}(\RR^s\UU^t)^{n-1}\RR^{s-1}$.  Such paths
are counted by Theorem~\ref{thm:periodic3}, above, in agreement with~\cite[Corollary 4]{chow}.

\newcommand{\conc}{\cdot}

\begin{figure}
    \includegraphics[width=.750\textwidth]{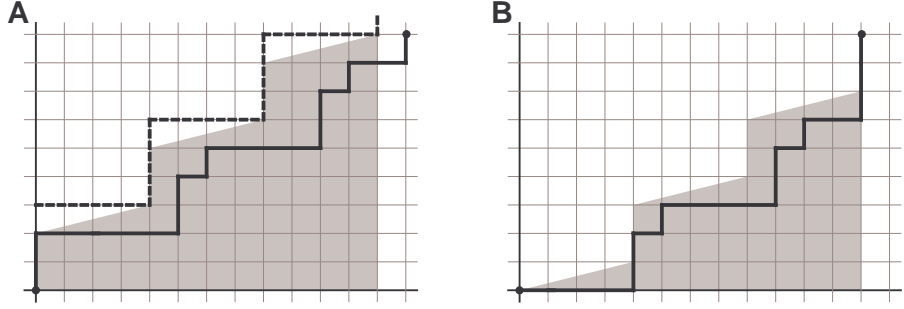}
    \caption{(A) The staircase $\stair{4}{3}$ (dotted line), a path $\pth{P}$ that avoids it (solid line), and
    the area dominated by the associated composition $\ba=(0^2,4)^3$ (shaded).  (B) The path $\pth{P}'$ is dominated by
    $\shift{\ba}{-2}=(4,0^2)^3$.
    }\label{fig:staircase2}
\end{figure}

Now consider a path $\pth{P}$ from $(0,0)$ to $(sn+1,tn)$ lying strictly below $\stair{s}{t}$. Clearly $\pth{P} =
\UU^j \RR \conc \pth{P}'$ for some $0 \leq j \leq t-1$ and some path $\pth{P}'$ from $(0,j)$ to $(sn,tn)$. Let
$\ba = (0^{t-1}, s)^n$.  Shift $\pth{P}'$ to the origin and append $j$ up steps to create the path $\pth{P''} =
\pth{P'} \conc \UU^j$ from $(0,0)$ to $(sn,tn)$.  (See Figure~\ref{fig:staircase2}.)
It is easy to check that $\pth{P}''$ is dominated by $\shift{\ba}{-j}$, and that every such path can be obtained
in this way. Thus there are $\sum_{j=0}^{t-1} D(\shift{\ba}{-j})$ paths to $(sn+1,tn)$ that avoid $\stair{s}{t}$.
From~\eqref{eq:periodic} and Theorem~\ref{thm:gen}, the sum evaluates to
$$
    \frac{1}{n}|\setgg{(sn,tn)}{\ba}| = \frac{1}{n}\binom{(s+t)n}{tn-1},
$$
again in accord with~\cite[Corollary 4]{chow}.

In fact,~\cite[Theorem 3]{chow} gives  formulae for the number of paths avoiding $\stair{s}{t}$ with a specified
number of corners.  For instance, performing the analysis above, but replacing Theorem~\ref{thm:gen} with the more
refined Theorem~\ref{thm:corners}, shows that the number of paths from $(0,0)$ to $(sn+1,tn)$ that avoid
$\stair{s}{t}$ and have $c$ up-right corners is
$$
    t\binom{sn}{c-1}\binom{tn}{c-1} - s\binom{sn-1}{c-2}\binom{tn+1}{c}.
$$
Note that we have used $c-1$ instead of $c$ in Theorem~\ref{thm:corners}, since the mapping $\pth{P} \mapsto
\pth{P}''$ described above reduces the number of up-right corners by 1.



\section*{Acknowledgements}

JI would like to thank Hugh Thomas for some helpful discussions, and in particular for suggesting a proof of
Corollary~\ref{cor:total} upon which the argument for Lemma~\ref{lem:cycleproof} is based. AR is supported by
a \emph{Natural Sciences and Engineering Research Council of Canada} Postdoctoral Fellowship.

\bibliographystyle{alpha}

\end{document}